\documentclass[11pt]{amsart}

\usepackage[a4paper,margin=3.2cm]{geometry}
\usepackage[T1]{fontenc}
\usepackage[utf8]{inputenc}
\usepackage{lmodern}
\usepackage{amsmath,amssymb,amsthm,mathtools}
\usepackage{mathrsfs}
\usepackage{hyperref}

\hypersetup{
  colorlinks=true,
  linkcolor=blue,
  citecolor=blue,
  urlcolor=blue
}

\newtheorem{theorem}{Theorem}[section]
\newtheorem{proposition}[theorem]{Proposition}
\newtheorem{lemma}[theorem]{Lemma}
\newtheorem{corollary}[theorem]{Corollary}
\theoremstyle{definition}

\theoremstyle{remark}
\newtheorem{remark}[theorem]{Remark}

\newcommand{\C}{\mathbb C}
\newcommand{\R}{\mathbb R}
\newcommand{\D}{\mathbb D}
\newcommand{\T}{\mathbb T}
\newcommand{\B}{\mathbb B}
\newcommand{\Sph}{\mathbb S}
\newcommand{\inner}[2]{\langle #1,#2\rangle}
\newcommand{\norm}[1]{\left\|#1\right\|}
\newcommand{\dist}{\operatorname{dist}}
\newcommand{\Aut}{\operatorname{Aut}}
\newcommand{\arcosh}{\operatorname{arcosh}}
\newcommand{\artanh}{\operatorname{artanh}}
\newcommand{\re}{\operatorname{Re}}

\title[Boundary Schwarz lemmas for holomorphic and minimal disks]
{Boundary Schwarz lemmas for holomorphic and conformal minimal disks}

\author{David Kalaj}
\address{University of Montenegro, Faculty of Natural Sciences and Mathematics, 81000 Podgorica, Montenegro}
\email{davidk@ucg.ac.me}

\date{September 11, 2025}

\begin{document}

\begin{abstract}
We prove two boundary Schwarz lemmas for disks. The first one concerns holomorphic maps from the unit disk into the unit ball of \(\C^m\). If \(F:\D\to\B_m\) is holomorphic, \(F(\zeta)\in\partial\B_m\) for some \(\zeta\in\T\), and a finite angular derivative exists at that point, then \(\norm{F'(\zeta)}\) admits an explicit lower bound in terms of \(F(0)\) and \(F'(0)\). The estimate is sharp, and we identify all equality maps; in particular, the extremals are one-dimensional Blaschke-type disks. The second result is a boundary Schwarz lemma for conformal minimal disks \(F:\D\to\B^n\subset\R^n\), \(n\ge3\). It follows from the distance-decreasing theorem of Forstneri\v c and Kalaj for the Poincare metric on the disk and the Cayley--Klein metric on the ball. We also determine the equality case: equality forces the image to be a totally geodesic planar disk, and in the noncentered case the boundary point is aligned with the radial direction of the center.
\end{abstract}

\maketitle

\section{Introduction}

The classical Schwarz lemma and the Schwarz--Pick lemma describe the metric contraction properties of holomorphic self-maps of the unit disk. Boundary forms of these results are equally important: they give quantitative information at boundary fixed points or, more generally, at boundary contact points. A sharp scalar boundary estimate was obtained by Osserman \cite{osserman2000}; related boundary Schwarz lemmas and rigidity results in several complex variables were developed by Krantz \cite{Stefan2011} and by Burns and Krantz \cite{bk1994}.

The first purpose of this paper is to prove a boundary Schwarz lemma for holomorphic disks in the Euclidean unit ball
\[
        \B_m=\{w\in\C^m:\norm{w}<1\}.
\]
The estimate is expressed not only in terms of the center \(F(0)\), but also in terms of the first derivative \(F'(0)\). This gives a vector-valued version of Osserman's boundary estimate and includes the sharp scalar equality case.

\begin{theorem}[Holomorphic boundary estimate and equality]\label{th:holomorphic-main}
Let \(F:\D\to\B_m\) be holomorphic. Assume that \(F\) has a nontangential boundary value \(F(\zeta)\in\partial\B_m\) at some \(\zeta\in\T\), and that \(F\) has a finite angular derivative \(F'(\zeta)\) there. Then
\begin{equation}\label{eq:main-holomorphic-bound}
        \norm{F'(\zeta)}
        \ge
        \frac{2(1-\norm{F(0)})^2}
        {1-\norm{F(0)}^2+\norm{F'(0)}}.
\end{equation}
The estimate is sharp.

More precisely, after precomposing with a rotation of the disk, we may assume that \(\zeta=1\). Equality in \eqref{eq:main-holomorphic-bound} holds if and only if there are a vector \(u\in\partial\B_m\), a number \(r\in[0,1)\), and a number \(t\in[0,1]\) such that
\[
        F(1)=u,
        \qquad
        F(0)=ru,
\]
and
\begin{equation}\label{eq:holomorphic-extremals}
        F(z)=u\,\frac{r+B_t(z)}{1+rB_t(z)},
        \qquad z\in\D,
\end{equation}
where
\[
        B_t(z)=z\frac{z+t}{1+tz},\qquad 0\le t<1,
        \quad\text{and}\quad B_1(z)=z.
\]
For these maps,
\[
        \norm{F'(0)}=(1-r^2)t
\]
and
\[
        \norm{F'(1)}
        =\frac{1-r}{1+r}\frac{2}{1+t}
        =\frac{2(1-r)^2}{1-r^2+\norm{F'(0)}}.
\]

In particular, if \(F(0)=0\), then
\begin{equation}\label{eq:centered-holomorphic-bound}
        \norm{F'(\zeta)}\ge \frac{2}{1+\norm{F'(0)}}.
\end{equation}
Hence \(\norm{F'(\zeta)}\ge1\) in the centered case. If \(F(0)=0\) and \(\norm{F'(\zeta)}=1\), then \(F\) is a complex linear disk:
\[
        F(z)=\bar\zeta z\,u,
        \qquad u\in\partial\B_m.
\]
\end{theorem}

The scalar equality case in the centered estimate is also rigid.

\begin{theorem}[Scalar equality case]\label{th:scalar-equality-intro}
Let \(f:\D\to\D\) be holomorphic and suppose that
\[
        f(0)=0,
        \qquad f(1)=1,
        \qquad
        |f'(1)|=\frac{2}{1+|f'(0)|}.
\]
Then either
\[
        f(z)=z,
\]
or
\[
        f(z)=z\frac{z+a}{1+az},
        \qquad 0\le a<1.
\]
Conversely, these functions give equality.
\end{theorem}

Osserman observed that the functions in Theorem~\ref{th:scalar-equality-intro} attain equality in the scalar boundary estimate \cite[Remark~2.2]{osserman2000}; this was also noted by Krantz \cite[Remark~5.2]{Stefan2011}. The converse is included here for completeness. A related version of the centered estimate was obtained by Ren and Wang \cite{RenWang2015}; see also Zhu \cite{Zhu2018} for boundary Schwarz lemmas for pluriharmonic mappings.

The second purpose of the paper is to prove the analogous boundary estimate for conformal minimal disks in real balls. Let
\[
        \B^n=\{x\in\R^n:\norm{x}<1\},
        \qquad
        \Sph^{n-1}=\partial\B^n.
\]
Forstneri\v c and Kalaj \cite{francdavid} proved that conformal minimal disks in \(\B^n\), \(n\ge3\), are distance decreasing from the Poincare metric on \(\D\) to the Cayley--Klein metric on \(\B^n\). We use this theorem to obtain the following boundary result.

\begin{theorem}[Boundary Schwarz lemma for conformal minimal disks]\label{th:minimal-main}
Let \(n\ge3\), and let \(F:\D\to\B^n\subset\R^n\) be a conformal minimal immersion. Assume that \(F\) has a boundary value \(F(\zeta)\in\Sph^{n-1}\) at some \(\zeta\in\T\), and that \(F\) extends differentiably to \(\zeta\) so that \(dF(\zeta)\) is defined. Then
\begin{equation}\label{eq:minimal-main-bound}
        \norm{dF(\zeta)}
        \ge
        \frac{1-\norm{F(0)}}{1+\norm{F(0)}}.
\end{equation}
If equality holds, then \(F(\D)\) is a totally geodesic linear disk. More precisely, there is a two-dimensional linear subspace \(L\subset\R^n\) such that
\[
        F(\D)=L\cap\B^n.
\]
In addition, if \(F(0)\ne0\), then
\[
        F(\zeta)=\frac{F(0)}{\norm{F(0)}}.
\]
Conversely, after identifying such an \(L\) with \(\C\), equality is attained, up to rotations of the source and of \(L\), by the disk automorphisms
\[
        z\longmapsto \frac{z+a}{1+az},
        \qquad 0\le a<1.
\]
\end{theorem}

For centered minimal disks this gives the following simple form.

\begin{corollary}[Centered case]\label{cor:centered-minimal-intro}
Let \(n\ge3\), and let \(F:\D\to\B^n\) be a conformal minimal immersion with \(F(0)=0\). If \(F\) extends differentiably to some \(\zeta\in\T\) and \(F(\zeta)\in\Sph^{n-1}\), then
\[
        \norm{dF(\zeta)}\ge1.
\]
Equality holds if and only if \(F\) parametrizes a totally geodesic planar disk through the origin.
\end{corollary}

The paper is organized as follows. Section~\ref{sec:holomorphic} proves the holomorphic boundary estimates, including the scalar equality theorem. Section~\ref{sec:minimal} proves the minimal-surface result and the equality case. The final corollary applies the boundary estimate to minimal disks in \(\R^3\) whose Gauss maps lie in a hemisphere.

\section{Holomorphic disks}\label{sec:holomorphic}

Throughout this section \(\B_m\) denotes the unit ball of \(\C^m\), and the Hermitian inner product is
\[
        \inner{x}{y}=\sum_{j=1}^m x_j\overline{y_j}.
\]

We begin with the centered estimate. It is the main input for the noncentered result.

\begin{theorem}[Centered holomorphic estimate]\label{th:centered-holomorphic}
Let \(F:\D\to\B_m\) be holomorphic and assume that \(F(0)=0\). Then
\begin{equation}\label{eq:interior-centered}
        \norm{F(z)}
        \le
        |z|\frac{|z|+\norm{F'(0)}}{1+|z|\norm{F'(0)}}
        \le |z|,
        \qquad z\in\D.
\end{equation}
If \(F\) has a nontangential boundary value \(F(\zeta)\in\partial\B_m\) at \(\zeta\in\T\), and if the finite angular derivative \(F'(\zeta)\) exists, then
\begin{equation}\label{eq:boundary-centered}
        \norm{F'(\zeta)}
        \ge
        \frac{2}{1+\norm{F'(0)}}.
\end{equation}
\end{theorem}

\begin{proof}
The vector-valued Schwarz lemma gives \(\norm{F(z)}\le |z|\) and \(\norm{F'(0)}\le1\). If \(\norm{F'(0)}=1\), the equality case in the Schwarz lemma implies \(F(z)=zv\) for some \(v\in\partial\B_m\), and the assertions follow. Hence assume that
\[
        A=\norm{F'(0)}<1.
\]
Define
\[
        f(z)=\frac{F(z)}{z},\qquad f(0)=F'(0).
\]
Since the linear extremal case has been excluded, the maximum principle gives \(f(\D)\subset\B_m\). Put \(a=f(0)\), \(w=f(z)\), \(x=\norm{w}\), and \(\rho=|z|\). By the Schwarz--Pick lemma for the ball, applied to the automorphism \(\varphi_a\in\Aut(\B_m)\) satisfying \(\varphi_a(a)=0\),
\[
        \norm{\varphi_a(w)}\le \rho.
\]
For the standard automorphism of the ball one has
\begin{equation}\label{eq:ball-identity}
        1-\norm{\varphi_a(w)}^2
        =
        \frac{(1-\norm{a}^2)(1-\norm{w}^2)}{|1-\inner{w}{a}|^2};
\end{equation}
see, for example, Rudin \cite[Ch.~2, Sec.~2.2]{Rudin2008}. Therefore
\[
        (1-A^2)(1-x^2)
        \ge
        (1-\rho^2)|1-\inner{w}{a}|^2
        \ge
        (1-\rho^2)(1-Ax)^2.
\]
Equivalently,
\[
        (1-A^2\rho^2)x^2-2A(1-\rho^2)x+(A^2-\rho^2)\le0.
\]
The roots of the corresponding quadratic polynomial are
\[
        \frac{A-\rho}{1-A\rho}
        \quad\text{and}\quad
        \frac{A+\rho}{1+A\rho}.
\]
It follows that
\[
        \norm{f(z)}\le \frac{A+|z|}{1+A|z|},
\]
and hence \eqref{eq:interior-centered} follows from \(F(z)=zf(z)\).

Assume now that \(F(\zeta)\in\partial\B_m\), \(|\zeta|=1\), and \(F'(\zeta)\) exists. For \(0<r<1\), \eqref{eq:interior-centered} gives
\[
        \norm{F(r\zeta)}
        \le
        r\frac{r+A}{1+Ar}.
\]
Since \(\norm{F(\zeta)}=1\),
\[
\begin{aligned}
        \left\|\frac{F(\zeta)-F(r\zeta)}{1-r}\right\|
        &\ge
        \frac{1-\norm{F(r\zeta)}}{1-r}  \\
        &\ge
        \frac{1-\dfrac{r(r+A)}{1+Ar}}{1-r}
        =
        \frac{1+r}{1+Ar}.
\end{aligned}
\]
Letting \(r\to1^-\) gives \eqref{eq:boundary-centered}.
\end{proof}

The following consequence is the centered rigidity statement in Theorem~\ref{th:holomorphic-main}.

\begin{corollary}\label{cor:centered-affine}
Let \(F:\D\to\B_m\) be holomorphic and \(F(0)=0\). Assume that \(F\) has a nontangential boundary value \(F(\zeta)\in\partial\B_m\) at \(\zeta\in\T\), and that the finite angular derivative \(F'(\zeta)\) exists. Then \(\norm{F'(\zeta)}\ge1\). If \(\norm{F'(\zeta)}=1\), then
\[
        F(z)=\bar\zeta z\,u
\]
for some \(u\in\partial\B_m\).
\end{corollary}

\begin{proof}
The inequality follows from Theorem~\ref{th:centered-holomorphic} and \(\norm{F'(0)}\le1\). If \(\norm{F'(\zeta)}=1\), then \eqref{eq:boundary-centered} implies \(\norm{F'(0)}=1\). Let \(v=F'(0)\). The scalar function \(g(z)=\inner{F(z)}{v}/\norm{v}\) maps \(\D\) into \(\D\), satisfies \(g(0)=0\), and has \(g'(0)=1\). The classical Schwarz lemma gives \(g(z)=z\). Equality in Cauchy's inequality then forces \(F(z)=zv\). The boundary condition at \(\zeta\) gives the displayed form.
\end{proof}

We next recall the standard ball automorphism and record the derivative estimates needed to remove the assumption \(F(0)=0\). Let \(a\in\B_m\), \(a\ne0\), and set
\[
        r=\norm{a},
        \qquad
        s=\sqrt{1-r^2},
        \qquad
        P_a w=\frac{\inner{w}{a}}{\norm{a}^2}a,
        \qquad
        Q_a=I-P_a.
\]
The standard involutive automorphism \(\varphi_a\) of \(\B_m\) satisfying \(\varphi_a(a)=0\) is
\begin{equation}\label{eq:ball-automorphism}
        \varphi_a(w)=\frac{a-P_aw-sQ_aw}{1-\inner{w}{a}}.
\end{equation}
For \(a=0\) we put \(\varphi_0(w)=-w\).

\begin{lemma}[Derivative estimates for ball automorphisms]\label{lem:automorphism-derivatives}
Let \(a\in\B_m\), \(r=\norm{a}\), and let \(\varphi_a\) be given by \eqref{eq:ball-automorphism}. Then
\begin{equation}\label{eq:derivative-at-a}
        \norm{D\varphi_a(a)}=\frac{1}{1-r^2},
\end{equation}
Moreover,
\begin{equation}\label{eq:derivative-at-zero}
        \norm{D\varphi_a(0)}\le\sqrt{1-r^2},
\end{equation}
with equality in \eqref{eq:derivative-at-zero} when \(m\ge2\); for \(m=1\) one has the sharper identity
\[
        \norm{D\varphi_a(0)}=1-r^2.
\]
and, for every \(w\in\overline{\B_m}\),
\begin{equation}\label{eq:global-derivative-bound}
        \norm{D\varphi_a(w)}\le \frac{1+r}{1-r}.
\end{equation}
\end{lemma}

\begin{proof}
The case \(a=0\) is immediate. Assume \(a\ne0\), put \(e=a/r\), and write every \(v\in\C^m\) as \(v=\alpha e+u\), where \(u\perp e\). Differentiating \eqref{eq:ball-automorphism} gives
\begin{equation}\label{eq:dphi-formula}
        D\varphi_a(w)[v]
        =
        \frac{-P_av-sQ_av+\inner{v}{a}\varphi_a(w)}{1-\inner{w}{a}}.
\end{equation}
At \(w=a\), since \(\varphi_a(a)=0\),
\[
        D\varphi_a(a)[v]
        =
        \frac{-P_av-sQ_av}{1-r^2},
\]
and the norm is \((1-r^2)^{-1}\). At \(w=0\), using \(\varphi_a(0)=a\), one obtains
\[
        D\varphi_a(0)[\alpha e+u]
        =-(1-r^2)\alpha e-su,
\]
Thus the norm is \(1-r^2\) when \(m=1\), while for \(m\ge2\) the transverse directions give the norm \(s=\sqrt{1-r^2}\). In particular, \eqref{eq:derivative-at-zero} holds in all dimensions. Finally, if \(\norm{w}\le1\) and \(\norm{v}=1\), then
\[
        \norm{-P_av-sQ_av}\le1,
        \qquad
        |\inner{v}{a}|\,\norm{\varphi_a(w)}\le r,
\]
while \(|1-\inner{w}{a}|\ge1-r\). Estimate \eqref{eq:global-derivative-bound} follows from \eqref{eq:dphi-formula}.
\end{proof}

The following exact formula is not needed for the main estimate, but it may be useful in other applications.

\begin{proposition}[Exact norm of \(D\varphi_a(w)\)]\label{prop:exact-automorphism-norm}
Let \(m\ge2\), \(a\in\B_m\setminus\{0\}\), \(r=\norm{a}\), \(s=\sqrt{1-r^2}\), and \(e=a/r\). Put \(y=\varphi_a(w)\),
\[
        b=-e+ry,
        \qquad b_\perp=(I-P_e)b,
\]
where \(P_e\) denotes orthogonal projection onto \(\C e\). Then
\begin{equation}\label{eq:exact-automorphism-norm}
\norm{D\varphi_a(w)}
=
\frac{1}{|1-\inner{w}{a}|}
\left(
\frac{
\norm{b}^2+s^2+\big((\norm{b}^2-s^2)^2+4s^2\norm{b_\perp}^2\big)^{1/2}
}{2}
\right)^{1/2}.
\end{equation}
In particular, \(\norm{D\varphi_a(a)}=(1-\norm{a}^2)^{-1}\).
\end{proposition}

\begin{proof}
Using \eqref{eq:dphi-formula} and writing \(v=\alpha e+u\), \(u\perp e\), we get
\[
        D\varphi_a(w)[v]
        =
        \frac{\alpha(-e+r\varphi_a(w))-su}{1-\inner{w}{a}}
        =
        \frac{\alpha b-su}{1-\inner{w}{a}}.
\]
Thus the problem is to compute the norm of
\[
        T(\alpha,u)=\alpha b-su,
        \qquad \alpha\in\C,\quad u\in e^\perp.
\]
Decompose \(b=b_\parallel+b_\perp\). Only the component of \(u\) in the direction of \(b_\perp\) can interact with \(\alpha b\); all components orthogonal to \(b_\perp\) contribute only \(s^2\norm{u}^2\). Hence \(\norm{T}^2\) is the largest eigenvalue of
\[
        \begin{pmatrix}
        \norm{b}^2 & -s\norm{b_\perp}\\
        -s\norm{b_\perp} & s^2
        \end{pmatrix}.
\]
The largest eigenvalue is the expression inside the outer square root in \eqref{eq:exact-automorphism-norm}. If \(w=a\), then \(y=0\), \(b=-e\), and \(b_\perp=0\), so \(\norm{T}=1\), giving the stated special case.
\end{proof}

We now prove the noncentered holomorphic estimate.

\begin{proof}[Proof of Theorem~\ref{th:holomorphic-main}]
By precomposing with a rotation of the disk, it is enough to prove the theorem for \(\zeta=1\). Put
\[
        a=F(0),
        \qquad
        r=\norm{a}.
\]
Let \(\varphi_a\) be the standard involutive automorphism of \(\B_m\) satisfying \(\varphi_a(a)=0\), and set
\[
        G=\varphi_a\circ F.
\]
Then \(G:\D\to\B_m\), \(G(0)=0\), and \(G(1)\in\partial\B_m\). The centered estimate applied to \(G\) gives
\[
        \norm{G'(1)}\ge \frac{2}{1+\norm{G'(0)}}.
\]
Since
\[
        G'(0)=D\varphi_a(a)F'(0)
\]
and
\[
        \norm{D\varphi_a(a)}=\frac1{1-r^2},
\]
we have
\[
        \norm{G'(0)}\le \frac{\norm{F'(0)}}{1-r^2}.
\]
Therefore
\begin{equation}\label{eq:main-proof-step-one}
        \norm{G'(1)}
        \ge
        \frac{2}{1+\dfrac{\norm{F'(0)}}{1-r^2}}
        =
        \frac{2(1-r^2)}{1-r^2+\norm{F'(0)}}.
\end{equation}
On the other hand,
\[
        G'(1)=D\varphi_a(F(1))F'(1).
\]
By \eqref{eq:global-derivative-bound},
\begin{equation}\label{eq:main-proof-step-two}
        \norm{G'(1)}
        \le
        \frac{1+r}{1-r}\,\norm{F'(1)}.
\end{equation}
Combining \eqref{eq:main-proof-step-one} and \eqref{eq:main-proof-step-two}, we obtain
\[
        \norm{F'(1)}
        \ge
        \frac{1-r}{1+r}\,
        \frac{2(1-r^2)}{1-r^2+\norm{F'(0)}}
        =
        \frac{2(1-r)^2}{1-r^2+\norm{F'(0)}}.
\]
This proves the inequality.

We now prove the equality statement. Assume that equality holds in \eqref{eq:main-holomorphic-bound}. Let
\[
        u=F(1)\in\partial\B_m.
\]
Then equality must hold in all the inequalities above. In particular, equality holds in the estimate
\[
        \norm{D\varphi_a(u)}\le \frac{1+r}{1-r}.
\]
Inspecting the proof of \eqref{eq:global-derivative-bound}, this forces
\[
        |1-\inner{u}{a}|=1-r.
\]
Since \(\norm{u}=1\) and \(\norm{a}=r\), this is possible only when
\[
        \inner{u}{a}=r.
\]
Thus equality holds in Cauchy's inequality, and hence
\begin{equation}\label{eq:equality-alignment-holomorphic}
        a=ru.
\end{equation}
So the center \(F(0)\) is aligned with the boundary value \(F(1)\).

Put
\[
        v=G(1)=\varphi_a(u).
\]
By \eqref{eq:equality-alignment-holomorphic}, the automorphism \(\varphi_a\) restricts to the complex line \(\C u\) as
\[
        \varphi_a(\lambda u)=u\frac{r-\lambda}{1-r\lambda}.
\]
In particular,
\begin{equation}\label{eq:v-minus-u}
        v=\varphi_a(u)=-u.
\end{equation}

Since equality holds for \(F\), equality also holds in the centered estimate for \(G\):
\begin{equation}\label{eq:centered-equality-for-G}
        \norm{G'(1)}=\frac{2}{1+\norm{G'(0)}}.
\end{equation}
Consider the scalar function
\[
        g(z)=\inner{G(z)}{v}.
\]
Then \(g:\D\to\D\), \(g(0)=0\), and \(g(1)=1\). By the Julia--Wolff--Caratheodory theorem, its angular derivative \(g'(1)\) is real and nonnegative. Moreover,
\[
        |g'(0)|\le \norm{G'(0)},
        \qquad
        |g'(1)|\le \norm{G'(1)}.
\]
The scalar centered boundary Schwarz lemma gives
\[
        |g'(1)|\ge \frac{2}{1+|g'(0)|}.
\]
Together with \eqref{eq:centered-equality-for-G}, we obtain
\[
        \frac{2}{1+\norm{G'(0)}}
        =\norm{G'(1)}
        \ge |g'(1)|
        \ge \frac{2}{1+|g'(0)|}
        \ge \frac{2}{1+\norm{G'(0)}}.
\]
Hence equality holds throughout. By Theorem~\ref{th:scalar-equality-intro},
\[
        g(z)=B_t(z),
        \qquad
        B_t(z)=z\frac{z+t}{1+tz},
        \qquad 0\le t<1,
\]
or \(g(z)=z\), which we denote by \(B_1(z)=z\). In this notation,
\[
        t=|g'(0)|=\norm{G'(0)}.
\]

It remains to show that \(G\) has no component orthogonal to \(v\). Write
\[
        G(z)=g(z)v+H(z),
        \qquad H(z)\perp v.
\]
Since \(G\) is bounded, its radial boundary values exist almost everywhere. The function \(g=B_t\) is a finite Blaschke product, and hence
\[
        |g(e^{i\theta})|=1
\]
for almost every \(e^{i\theta}\in\T\). Therefore, for almost every boundary point,
\[
        1\ge \norm{G(e^{i\theta})}^2
        =|g(e^{i\theta})|^2+\norm{H(e^{i\theta})}^2
        =1+\norm{H(e^{i\theta})}^2.
\]
It follows that \(H(e^{i\theta})=0\) almost everywhere, and hence \(H\equiv0\). Thus
\[
        G(z)=vB_t(z)=-uB_t(z),
\]
where the last equality follows from \eqref{eq:v-minus-u}.

Finally, using that \(\varphi_a\) is an involution and that \(a=ru\), we get
\[
        F(z)=\varphi_a(G(z))
        =\varphi_{ru}(-uB_t(z))
        =u\frac{r+B_t(z)}{1+rB_t(z)}.
\]
This proves the necessity of the extremal form \eqref{eq:holomorphic-extremals}.

Conversely, let \(F\) be of the form
\[
        F(z)=u\frac{r+B_t(z)}{1+rB_t(z)},
        \qquad
        u\in\partial\B_m,
        \quad 0\le r<1,
        \quad 0\le t\le1.
\]
Then \(F(0)=ru\), and a direct calculation gives
\[
        \norm{F'(0)}=(1-r^2)t
\]
and
\[
        \norm{F'(1)}
        =\frac{1-r}{1+r}\frac{2}{1+t}.
\]
Hence
\[
        \norm{F'(1)}
        =
        \frac{2(1-r)^2}{1-r^2+\norm{F'(0)}}.
\]
Thus equality is attained exactly by the maps in \eqref{eq:holomorphic-extremals}.

The centered conclusions follow by putting \(r=0\) and then restoring the original boundary point \(\zeta\). If \(F(0)=0\), the estimate becomes
\[
        \norm{F'(\zeta)}\ge \frac{2}{1+\norm{F'(0)}}.
\]
If, in addition, \(\norm{F'(\zeta)}=1\), then necessarily \(\norm{F'(0)}=1\). The equality case in the vector Schwarz lemma gives
\[
        F(z)=\bar\zeta z\,u,
        \qquad u\in\partial\B_m.
\]
This completes the proof.
\end{proof}

\begin{remark}\label{rem:explicit-noncentered}
The estimate \(\norm{F'(0)}\le\sqrt{1-\norm{F(0)}^2}\) follows by applying the ball Schwarz--Pick lemma to \(\varphi_{F(0)}\circ F\), together with the inequality \eqref{eq:derivative-at-zero}. In dimension \(m=1\) this step gives the sharper classical bound \(|F'(0)|\le1-|F(0)|^2\). Hence Theorem~\ref{th:holomorphic-main} implies the explicit lower bound
\[
        \norm{F'(\zeta)}
        \ge
        \frac{2(1-r)^2}{1-r^2+\sqrt{1-r^2}},
        \qquad r=\norm{F(0)},
\]
for maps into \(\B_m\). In the scalar case the classical Schwarz--Pick lemma gives the sharper estimate \(|F'(0)|\le1-|F(0)|^2\), and therefore
\[
        |F'(\zeta)|\ge\frac{1-|F(0)|}{1+|F(0)|}.
\]
\end{remark}

It remains to prove the scalar equality theorem. We use the Julia inequality in a standard form.

\begin{proposition}[Julia inequality]\label{prop:julia}
Let \(f:\D\to\D\) be holomorphic. Assume that \(f\) has nontangential limit \(1\) at \(1\) and a finite angular derivative \(f'(1)\) there. Then
\begin{equation}\label{eq:julia}
        \frac{|1-f(z)|^2}{1-|f(z)|^2}
        \le
        f'(1)\frac{|1-z|^2}{1-|z|^2},
        \qquad z\in\D.
\end{equation}
In particular,
\[
        f'(1)\ge \frac{|1-f(0)|^2}{1-|f(0)|^2}.
\]
If equality holds in \eqref{eq:julia} for one point \(z\in\D\), then \(f\) is a disk automorphism or a constant boundary map.
\end{proposition}

\begin{proof}
This is the classical Julia--Wolff--Caratheodory theorem; see Caratheodory \cite[p.~27]{JuliaCar}.
\end{proof}

\begin{proof}[Proof of Theorem~\ref{th:scalar-equality-intro}]
Write the zero of \(f\) at the origin with multiplicity \(m\ge1\):
\[
        f(z)=z^m h(z),
        \qquad h(1)=1.
\]
By the generalized Schwarz lemma, \(|h(z)|\le1\) in \(\D\). First note that the Julia--Wolff--Caratheodory theorem gives \(h'(1)\ge0\), and hence \(f'(1)=m+h'(1)\ge m\). Since the equality assumption gives \(f'(1)=2/(1+|f'(0)|)\le2\), necessarily \(m\le2\). If \(h\) is constant of modulus one, then \(h\equiv1\), so \(f(z)=z\) if \(m=1\) and \(f(z)=z^2\) if \(m=2\); the latter is the second family with \(a=0\). We may therefore assume that \(h:\D\to\D\).

The Julia--Wolff--Caratheodory theorem gives \(h'(1)\ge0\). Hence
\[
        f'(1)=m+h'(1)\ge m.
\]

If \(m=2\), then \(f'(0)=0\), so \(f'(1)=2\). Thus \(h'(1)=0\). Julia's inequality forces \(h\equiv1\), and hence \(f(z)=z^2\).

It remains to consider \(m=1\). Write \(f(z)=zh(z)\), and put \(a=h(0)=f'(0)\). If \(|a|=1\), the maximum principle gives \(h\equiv a\); since \(h(1)=1\), this gives \(a=1\) and \(f(z)=z\). Assume now that \(|a|<1\). From the equality assumption,
\[
        1+h'(1)=\frac{2}{1+|a|},
        \qquad\text{so}\qquad
        h'(1)=\frac{1-|a|}{1+|a|}.
\]
Julia's inequality applied to \(h\) gives
\[
        h'(1)
        \ge
        \frac{|1-a|^2}{1-|a|^2}.
\]
Writing \(a=\rho e^{it}\), comparison with the preceding identity implies \(\cos t=1\). Hence \(a\in[0,1)\). Equality holds in Julia's inequality, so \(h\) is the automorphism of \(\D\) fixing \(1\) and satisfying \(h(0)=a\), namely
\[
        h(z)=\frac{z+a}{1+az},
        \qquad 0\le a<1.
\]
Thus
\[
        f(z)=z\frac{z+a}{1+az},
        \qquad 0\le a<1.
\]
This proves the theorem.
\end{proof}

\section{Conformal minimal disks}\label{sec:minimal}

Let \(\B^n\subset\R^n\) be the Euclidean unit ball. The Cayley--Klein distance on \(\B^n\) is
\begin{equation}\label{eq:ck-distance}
        \dist(x,y)
        =
        \arcosh\left(
        \frac{|1-\inner{x}{y}|}{\sqrt{(1-\norm{x}^2)(1-\norm{y}^2)}}
        \right),
        \qquad x,y\in\B^n.
\end{equation}
We shall use the following theorem of Forstneri\v c and Kalaj.

\begin{theorem}[Distance-decreasing property]\label{th:distance-decreasing}
Let \(n\ge3\), and let \(F:\D\to\B^n\) be a conformal minimal immersion. Then
\begin{equation}\label{eq:distance-decreasing}
        \dist(F(z),F(w))\le \dist_{\mathcal P}(z,w),
        \qquad z,w\in\D,
\end{equation}
where \(\dist_{\mathcal P}\) is the Poincare distance on \(\D\). If equality holds for one pair of distinct points, then \(F(\D)\) is a totally geodesic affine disk, that is, the intersection of \(\B^n\) with a two-dimensional affine plane; in this case equality holds for all pairs of points.
\end{theorem}

\begin{proof}
This is \cite[Theorem~1.4]{francdavid}; see also \cite{DrinovecDrnovsekForstneric2023} for the associated minimal metric.
\end{proof}

We first extract a Euclidean consequence of \eqref{eq:distance-decreasing}.

\begin{lemma}\label{lem:minimal-growth}
Let \(n\ge3\), and let \(F:\D\to\B^n\) be a conformal minimal immersion. Then
\begin{equation}\label{eq:minimal-growth}
        \norm{F(z)}
        \le
        \frac{|z|+\norm{F(0)}}{1+|z|\norm{F(0)}}.
\end{equation}
In particular, if \(F(0)=0\), then \(\norm{F(z)}\le |z|\).
\end{lemma}

\begin{proof}
Let \(p=F(0)\). By Theorem~\ref{th:distance-decreasing},
\[
        \dist(F(z),p)
        \le
        \dist_{\mathcal P}(z,0)
        =\artanh |z|.
\]
We spell out the last step. Put \(a=\norm p\), \(\rho=|z|\), and \(y=F(z)\). If \(a=0\), the assertion follows directly. Assume \(a>0\) and rotate so that \(p=ae_1\). Since
\[
        \cosh(\dist(y,p))
        =
        \frac{1-\inner{p}{y}}
        {\sqrt{(1-a^2)(1-\norm y^2)}}
        \le
        \cosh(\artanh\rho)=\frac1{\sqrt{1-\rho^2}},
\]
we get
\[
        (1-\rho^2)(1-a y_1)^2
        \le (1-a^2)(1-\norm y^2).
\]
For fixed \(R=\norm y\), the left-hand side is minimized when \(y_1=R\). Hence a necessary condition for the existence of such \(y\) is
\[
        (1-\rho^2)(1-aR)^2
        \le (1-a^2)(1-R^2).
\]
Solving this one-dimensional inequality gives
\[
        R\le \frac{a+\rho}{1+a\rho}.
\]
Thus \(\norm{F(z)}=R\le(a+\rho)/(1+a\rho)\), which is \eqref{eq:minimal-growth}.
\end{proof}

We now prove the minimal boundary Schwarz lemma.

\begin{proof}[Proof of Theorem~\ref{th:minimal-main}]
After precomposing with a rotation of \(\D\), we may assume that \(\zeta=1\). Put
\[
        p=F(0),
        \qquad
        a=\norm{p},
        \qquad
        c=\frac{1-a}{1+a},
        \qquad
        q=F(1)\in\Sph^{n-1}.
\]
For \(0<r<1\), write
\[
        y_r=F(r),
        \qquad
        \rho(r)=\norm{y_r}.
\]
By Lemma~\ref{lem:minimal-growth},
\begin{equation}\label{eq:rho-bound}
        \rho(r)\le \frac{r+a}{1+ar}.
\end{equation}
Since \(\norm{q}=1\),
\[
        \norm{q-y_r}\ge 1-\rho(r).
\]
Therefore
\[
\begin{aligned}
        \left\|\frac{F(1)-F(r)}{1-r}\right\|
        &\ge
        \frac{1-\rho(r)}{1-r} \\
        &\ge
        \frac{1-\dfrac{r+a}{1+ar}}{1-r}
        =
        \frac{1-a}{1+ar}.
\end{aligned}
\]
Letting \(r\to1^-\) and using the existence of \(dF(1)\) gives
\[
        \norm{dF(1)[1]}\ge c.
\]
Since \(F\) is conformal, the norm of \(dF(1)\) is independent of the unit tangent direction in the source. Hence
\[
        \norm{dF(1)}\ge c,
\]
which proves \eqref{eq:minimal-main-bound}.

Assume now that equality holds. Differentiability at \(1\) and conformality give
\begin{equation}\label{eq:boundary-asymptotic}
        \frac{\norm{F(1)-F(r)}}{1-r}\longrightarrow c,
        \qquad r\to1^-.
\end{equation}
The preceding chain of inequalities and \((1-a)/(1+ar)\to c\) imply
\begin{equation}\label{eq:rho-asymptotic}
        \frac{1-\rho(r)}{1-r}\longrightarrow c.
\end{equation}

We compare \(\dist(p,y_r)\) with \(\dist_{\mathcal P}(0,r)\). By \eqref{eq:ck-distance},
\[
        \dist(p,y_r)
        =
        \arcosh\left(
        \frac{1-\inner{p}{y_r}}
        {\sqrt{(1-a^2)(1-\rho(r)^2)}}
        \right),
\]
while
\[
        \dist_{\mathcal P}(0,r)=\frac12\log\frac{1+r}{1-r}.
\]
Since \(y_r\to q\in\Sph^{n-1}\), \(\rho(r)\to1\), and \(1-\inner{p}{q}>0\), the argument of \(\arcosh\) tends to \(+\infty\). Using \(\arcosh X=\log(2X)+o(1)\), \(X\to+\infty\), and \eqref{eq:rho-asymptotic}, we obtain
\begin{equation}\label{eq:distance-defect-limit}
        \lim_{r\to1^-}
        \bigl(\dist_{\mathcal P}(0,r)-\dist(p,F(r))\bigr)
        =
        \log\frac{1-a}{1-\inner{p}{q}}.
\end{equation}
The distance-decreasing theorem gives \(\dist(p,F(r))\le\dist_{\mathcal P}(0,r)\), so the limit in \eqref{eq:distance-defect-limit} is nonnegative. Hence
\[
        \inner{p}{q}\ge a.
\]
By Cauchy's inequality,
\[
        \inner{p}{q}\le \norm{p}\,\norm{q}=a.
\]
Thus
\begin{equation}\label{eq:radial-alignment}
        \inner{p}{q}=a.
\end{equation}
If \(a>0\), equality in Cauchy's inequality gives
\[
        q=\frac{p}{\norm{p}},
\]
which is the asserted radial alignment. Equations \eqref{eq:distance-defect-limit} and \eqref{eq:radial-alignment} also give
\begin{equation}\label{eq:distance-defect-zero}
        \dist_{\mathcal P}(0,r)-\dist(p,F(r))\longrightarrow0,
        \qquad r\to1^-.
\end{equation}

Fix \(s\in(0,1)\). For \(s<r<1\), the triangle inequality and Theorem~\ref{th:distance-decreasing} give
\[
\begin{aligned}
        \dist(p,F(r))
        &\le \dist(p,F(s))+\dist(F(s),F(r)) \\
        &\le \dist(p,F(s))+\dist_{\mathcal P}(s,r).
\end{aligned}
\]
Since \(0<s<r<1\) lie on a Poincare geodesic,
\[
        \dist_{\mathcal P}(0,r)
        =\dist_{\mathcal P}(0,s)+\dist_{\mathcal P}(s,r).
\]
Consequently,
\[
        \dist_{\mathcal P}(0,r)-\dist(p,F(r))
        \ge
        \dist_{\mathcal P}(0,s)-\dist(p,F(s)).
\]
Letting \(r\to1^-\) and using \eqref{eq:distance-defect-zero}, we get
\[
        \dist_{\mathcal P}(0,s)
        \le
        \dist(p,F(s)).
\]
The reverse inequality is Theorem~\ref{th:distance-decreasing}. Therefore
\[
        \dist(p,F(s))=\dist_{\mathcal P}(0,s).
\]
Thus equality holds in Theorem~\ref{th:distance-decreasing} for the distinct pair \(0,s\). Hence \(F(\D)\) is a totally geodesic affine disk,
\[
        F(\D)=A\cap\B^n,
\]
where \(A\subset\R^n\) is a two-dimensional affine plane.

It remains to see that \(A\) passes through the origin. If \(a=0\), then \(p=0\in A\). If \(a>0\), then both \(p\) and \(q=p/\norm{p}\) belong to \(A\), and the line through \(p\) and \(q\) contains the origin. Hence \(0\in A\). Therefore \(A=L\) is a two-dimensional linear subspace of \(\R^n\), and
\[
        F(\D)=L\cap\B^n.
\]

Conversely, identify \(L\) with \(\C\). Up to rotations of the source and of \(L\), the disk automorphism
\[
        z\longmapsto \frac{z+a}{1+az},
        \qquad 0\le a<1,
\]
satisfies \(F(0)=aF(1)\) and
\[
        \norm{dF(1)}=\frac{1-a}{1+a}
        =
        \frac{1-\norm{F(0)}}{1+\norm{F(0)}}.
\]
This proves the equality statement and completes the proof.
\end{proof}

\begin{proof}[Proof of Corollary~\ref{cor:centered-minimal-intro}]
If \(F(0)=0\), Theorem~\ref{th:minimal-main} gives \(\norm{dF(\zeta)}\ge1\). Its equality statement is exactly the asserted characterization; the radial alignment condition is void in the centered case.
\end{proof}

We finish with a consequence for conformal minimal disks in \(\R^3\) whose normal vectors lie in a hemisphere. After a rigid motion, such a disk has Enneper--Weierstrass data \((p,q)\) with \(|q|<1\):
\begin{equation}\label{eq:weierstrass-data}
        F(z)=\re\int^z
        \left(
        \frac12(1-q^2)p,
        \frac{i}{2}(1+q^2)p,
        qp
        \right).
\end{equation}
Here \(q\) is the stereographic representation of the Gauss map and \(p\) is a holomorphic one-form. With stereographic projection taken from the south pole, the unit normal is
\[
        \mathbf n(z)=\frac{1}{1+|q(z)|^2}
        \bigl(2\re q(z),2\operatorname{Im}q(z),1-|q(z)|^2\bigr),
\]
so \(n_3>0\) is equivalent to \(|q|<1\). For \(z=re^{it}\), with the normalization in \eqref{eq:weierstrass-data}, one has
\begin{equation}\label{eq:radial-factor}
        |F_r(z)|=\norm{dF(z)}=\frac12 |p(z)|(1+|q(z)|^2),
\end{equation}
see, for example, Duren \cite[Ch.~9]{Duren2004}.

\begin{corollary}\label{cor:hemisphere-lipschitz}
Let \(F:\D\to\B^3\) be a conformal minimal embedding whose unit normals belong to an open hemisphere. Assume that \(F\) extends differentiably to \(\T\), that \(F(\T)\subset\Sph^2\), and that the Weierstrass data in \eqref{eq:weierstrass-data} extend continuously to \(\overline\D\). Then
\begin{equation}\label{eq:interior-lower-bound}
        \norm{dF(z)}
        \ge
        \frac12\frac{1-\norm{F(0)}}{1+\norm{F(0)}},
        \qquad z\in\D.
\end{equation}
Consequently, \(F^{-1}:F(\D)\to\D\) is Lipschitz continuous with respect to the intrinsic metric on \(F(\D)\) and the Euclidean metric on \(\D\).
\end{corollary}

\begin{proof}
Theorem~\ref{th:minimal-main} gives, for every \(\zeta\in\T\),
\[
        \norm{dF(\zeta)}
        \ge
        \frac{1-\norm{F(0)}}{1+\norm{F(0)}}.
\]
Since \(F\) is conformal, this norm is the conformal factor. In particular, by \eqref{eq:radial-factor},
\[
        \norm{dF(z)}=|F_r(z)|=\frac12 |p(z)|(1+|q(z)|^2).
\]
Since \(|q|<1\), we have \(\frac12(1+|q|^2)\le1\). The boundary estimate therefore gives
\[
        |p(\zeta)|
        \ge
        \frac{1-\norm{F(0)}}{1+\norm{F(0)}}.
\]
The immersion condition gives \(p\ne0\) in \(\D\). Applying the maximum principle to \(1/p\) yields the same lower bound for \(|p|\) throughout \(\D\). Hence \eqref{eq:interior-lower-bound} follows from \eqref{eq:radial-factor}, because \(\frac12(1+|q|^2)\ge\frac12\).

If \(\gamma\) is a \(C^1\) curve in \(\D\), then
\[
        \operatorname{length}_{F(\D)}(F\circ\gamma)
        \ge
        \frac12\frac{1-\norm{F(0)}}{1+\norm{F(0)}}
        \operatorname{length}_{\D}(\gamma).
\]
Taking infima over curves joining two points gives
\[
        |F^{-1}(x)-F^{-1}(y)|
        \le
        2\frac{1+\norm{F(0)}}{1-\norm{F(0)}}
        d_{F(\D)}(x,y),
\]
where \(d_{F(\D)}\) is the intrinsic distance on the minimal disk. This proves the Lipschitz assertion.
\end{proof}

\subsection*{Funding}
The author is partially supported by the Ministry of Education, Science and Innovation of Montenegro through the grants \emph{Mathematical Analysis, Optimization and Machine Learning} and \emph{Complex-analytic and geometric techniques for non-Euclidean machine learning: theory and applications}.

\subsection*{Data availability}
Data sharing is not applicable to this article since no data sets were generated or analyzed.

\subsection*{Ethics declarations}
The author declares that he has no conflict of interest.


\begin{thebibliography}{99}

\bibitem{abate}
M. Abate and J. Raissy,
\newblock A Julia--Wolff--Caratheodory theorem for infinitesimal generators in the unit ball,
\newblock \emph{Trans. Amer. Math. Soc.} \textbf{368} (2016), no. 8, 5415--5431.

\bibitem{bk1994}
D. M. Burns and S. G. Krantz,
\newblock Rigidity of holomorphic mappings and a new Schwarz lemma at the boundary,
\newblock \emph{J. Amer. Math. Soc.} \textbf{7} (1994), no. 3, 661--676.

\bibitem{JuliaCar}
C. Caratheodory,
\newblock \emph{Conformal Representation}, 2nd ed.,
\newblock Cambridge Tracts in Mathematics and Mathematical Physics, Vol. 28,
\newblock Cambridge University Press, Cambridge, 1952.

\bibitem{Cannon1997}
J. W. Cannon, W. J. Floyd, R. Kenyon and W. R. Parry,
\newblock Hyperbolic geometry,
\newblock in \emph{Flavors of Geometry}, Math. Sci. Res. Inst. Publ., Vol. 31, pp. 59--115,
\newblock Cambridge University Press, Cambridge, 1997.

\bibitem{DrinovecDrnovsekForstneric2023}
B. Drinovec Drnov\v sek and F. Forstneri\v c,
\newblock Hyperbolic domains in real Euclidean spaces,
\newblock \emph{Pure Appl. Math. Q.} \textbf{19} (2023), no. 6, 2689--2735.

\bibitem{Duren2004}
P. Duren,
\newblock \emph{Harmonic Mappings in the Plane},
\newblock Cambridge Tracts in Mathematics, Vol. 156,
\newblock Cambridge University Press, Cambridge, 2004.

\bibitem{francdavid}
F. Forstneri\v c and D. Kalaj,
\newblock Schwarz--Pick lemma for harmonic maps which are conformal at a point,
\newblock \emph{Anal. PDE} \textbf{17} (2024), no. 3, 981--1003.

\bibitem{glasgow}
D. Kalaj,
\newblock Schwarz lemma for holomorphic mappings in the unit ball,
\newblock \emph{Glasg. Math. J.} \textbf{60} (2018), no. 1, 219--224.

\bibitem{Stefan2011}
S. G. Krantz,
\newblock The Schwarz lemma at the boundary,
\newblock \emph{Complex Var. Elliptic Equ.} \textbf{56} (2011), no. 5, 455--468.

\bibitem{osserman2000}
R. Osserman,
\newblock A sharp Schwarz inequality on the boundary,
\newblock \emph{Proc. Amer. Math. Soc.} \textbf{128} (2000), no. 12, 3513--3517.

\bibitem{RenWang2015}
G. Ren and X. Wang,
\newblock Extremal functions of boundary Schwarz lemma,
\newblock arXiv:1502.02369, 2015.

\bibitem{Rudin2008}
W. Rudin,
\newblock \emph{Function Theory in the Unit Ball of \(\C^n\)},
\newblock Classics in Mathematics,
\newblock Springer-Verlag, Berlin, 2008.

\bibitem{Zhu2018}
J.-F. Zhu,
\newblock Schwarz lemma and boundary Schwarz lemma for pluriharmonic mappings,
\newblock \emph{Filomat} \textbf{32} (2018), no. 15, 5385--5402.

\end{thebibliography}
\end{document}